\documentclass[11pt]{article}
\usepackage{amsthm}
\usepackage{latexsym, amssymb, amsmath}
\usepackage{graphicx}
\usepackage[dvipsnames]{xcolor}

\begin{document}
\numberwithin{equation}{section}

%%%%%%%%%%%%%%%%%%%%%%%%%%%%%%%%%%%%%%%%%%%%%%%%%%%%%%%%%%%%%%%%
\newtheorem{THEOREM}{Theorem}
\newtheorem{PRO}{Proposition}
\newtheorem{XXXX}{\underline{Theorem}}
\newtheorem{CLAIM}{Claim}
\newtheorem{COR}{Corollary}
\newtheorem{LEMMA}{Lemma}
\newtheorem{REM}{Remark}
\newtheorem{EX}{Example}
\newenvironment{PROOF}{{\bf Proof}.}{{\ \vrule height7pt width4pt depth1pt} \par \vspace{2ex} }
\newcommand{\Bibitem}[1]{\bibitem{#1} \ifnum\thelabelflag=1 
  \marginpar{\vspace{0.6\baselineskip}\hspace{-1.08\textwidth}\fbox{\rm#1}}
  \fi}
\newcounter{labelflag} \setcounter{labelflag}{0}
\newcommand{\labelon}{\setcounter{labelflag}{1}}
\newcommand{\Label}[1]{\label{#1} \ifnum\thelabelflag=1 
  \ifmmode  \makebox[0in][l]{\qquad\fbox{\rm#1}}
  \else\marginpar{\vspace{0.7\baselineskip}\hspace{-1.15\textwidth}\fbox{\rm#1}}
  \fi \fi}
% \labelon

\newcommand{\LEFTLINE}{\ifhmode\newline\else\noindent\fi}
\newcommand{\RIGHTLINE}[1]{\LEFTLINE\rightline{#1}}
\newcommand{\CENTERLINE}[1]{\LEFTLINE\centerline{#1}}
\def\BOX #1 #2 {\framebox[#1in]{\parbox{#1in}{\vspace{#2in}}}}
\parskip=8pt plus 2pt
\def\AUTHOR#1{\author{#1} \maketitle}
\def\Title#1{\begin{center}  \Large\bf #1 \end{center}  \vskip 1ex }
\def\Author#1{\vspace*{-2ex}\begin{center} #1 \end{center}  
 \vskip 2ex \par}
\renewcommand{\theequation}{\arabic{section}.\arabic{equation}}
\def\bdk#1{\makebox[0pt][l]{#1}\hspace*{0.03ex}\makebox[0pt][l]{#1}\hspace*{0.03ex}\makebox[0pt][l]{#1}\hspace*{0.03ex}\makebox[0pt][l]{#1}\mbox{#1} }
\def\psbx#1 #2 {\mbox{\psfig{file=#1,height=#2}}}

%pm  /po  polynomial
%pm  /py  positivity
%pm  /PS  partial sum
%pm  /ni  non-increasing
\newcommand{\FG}[2]{{\includegraphics[height=#1mm]{#2.eps}}}
 
%%%%%%%%%%%%%%%%%%%%%%%%%%%%%%%%%%%%%%%%%%%%%%%%%%%%%%%%%%%%%%%%
 
\Title{Improved Vietoris Sine Inequalities for \\  Non-Monotone, Non-Decaying Coefficients}

\par\vspace*{-4mm}\par
\Author{(Draft version 1, April 25, 2015)}
% \vspace{0.5cm}
\begin{center}
MAN KAM KWONG\footnote{The research of this author is supported by the Hong Kong Government GRF Grant PolyU 5003/12P and the Hong Kong Polytechnic University Grants G-UC22 and G-UA10}
\end{center}

% \vspace{0.5cm}
\begin{center}
\emph{Department of Applied Mathematics\\ The Hong Kong Polytechnic University, Hunghom, Hong Kong}\\
\tt{mankwong@polyu.edu.hk}
\end{center}

\par\vspace*{\baselineskip}\par

\newcommand{\mb}{\mathbf}

\begin{abstract}
\parskip=6pt
The classical Vietoris sine inequality states that for any non-increasing
sequence of positive real numbers $ \left\{ a_k\right\} _{k=1}^\infty  $ satisfying
$$  \hspace*{20mm}   a_{2j-1} \, \geq  \frac{2j}{2j-1} \,\,a_{2j} \qquad  (j=1,2,3,\cdots) ,  \eqno (*)  $$
the following sine polynomials are nonnegative in $ [0,\pi ] $,
$$  \hspace*{9mm}  \sum_{k=1}^{n} \, a_k \, \sin(kx) \geq  0, \qquad  x\in[0,\pi ], \quad \mbox{for all } n = 1,2,3,\cdots .  \eqno (\dagger)  $$

Recently, the author has improved this result to include non-monotone sequences.

In this paper, we establish two further extensions.
The first states that
if $ \left\{ a_k\right\}  $ is a sequence of positive numbers satisfying
$$  \hspace*{20mm}     a_{2}\, \geq \, 0.5869890995\cdots  \,\, a_{3}, \quad  \mbox{and} \quad  a_{2j}\, \geq  \,  \frac{2j+1}{2j+2} \,\,a_{2j+1} \qquad  (j=2,3,\cdots) ,  $$
then ($*$) implies ($\dagger$).
An example is
$ \left\{ a_k\right\}  =          \left\{  \frac{8}{5} , \frac{4}{5} , \frac{4}{3} , 1, \frac{6}{5} , 1, \frac{8}{7} , 1,  \cdots \right\}  , $
with $ a_k=1 $ for even $ k\geq 4 $ and $ a_k=(k+1)/k $ for odd $ k\geq 3 $.

A second, independent, extension affirms that ($\dagger$) also holds under
($*$) and
$$  \hspace*{20mm}    a_{2j}\, \geq  \, \frac{(2j+1)(4j-1)}{2j(4j+3)} \,\,a_{2j+1} \qquad  (j=1,2,\cdots).  $$
An example is
$ \left\{  3, \frac{3}{2} , \frac{7}{3} , \frac{7}{4} ,  \frac{11}{5} , \frac{11}{6} ,  \cdots \right\}   $
where $ a_{k}=2-\frac{(-1)^k}{k} $.

The coefficients in these examples are not monotone and not converging to 0.
\end{abstract}

\vspace{1cm}


\newpage

\section{Introduction}

Excellent surveys on the history and applications
of nonnegative trigonometric polynomials can be found, for example, in Alzer, Koumandos
and Lamprecht \cite{akl},
Askey et.\ al.\ \cite{A}--\cite{AS}, Brown \cite{b1}, and Koumandos \cite{Ks2},
and the references therein.

For convenience, we use the acronyms NN to stand for ``non-negative'', and
PS for P-Sum (a sum with all its partial sums NN). These can be 
interpreted as an adjective or a noun depending on the context.
A sequence of real numbers is denoted by 
$ \left\{ a_k\right\} _{k=1}^\infty  $, or simply, $ \left\{ a_k\right\}  $. A finite $ n $-tuple of numbers can
be interpreted as an infinite sequence by adding $ 0 $ to the end.
The symbol $ \searrow $ means non-increasing. 

Following the convention adopted in \cite{Kw2}, we use bold capital letters
such as $ \mb F $ and $ \mb\Phi  $ to denote sums of numbers or functions.
One of the earliest known PS is
\begin{equation}  \mb F = \sum \frac{\sin(kx)}{k} \,,  \Label{F}  \end{equation}
first conjectured by Ferj\'er 1910, and confirmed 
independently by Jackson and Gronwall. 
Vietoris, in 1958, established a deep result that includes $ \mb F $.

{
\renewcommand{\theTHEOREM}{A}
\begin{THEOREM}[Vietoris \cite{V}]
The sum \,
$ \sum a_k\,\sin(kx) $ is a PS in $ [0,\pi ] $ $($i.e. $(\dagger)$ holds$)$ if
$a_k \searrow 0$ and
\begin{equation}  a_{2j-1} \,\, \geq  \,\, \frac{2j}{2j-1} \,\, a_{2j} \,, \qquad  \mbox{for} \quad  j=1,2,\cdots \,.  \Label{v}  \end{equation}
\end{THEOREM}
}

\begin{REM} \em
There is an analogous cosine inequality (%
$ a_1 + \sum a_k\,\cos(kx)  $
is also a PS), but we are only concerned with
the sine sum in this paper.
\end{REM}

Belov, in 1995, greatly improved Vietoris' sine inequality, 
by establishing, under the monotonicity requirement, a necessary and sufficient 
condition for PS.

{
\renewcommand{\theTHEOREM}{B}
\begin{THEOREM}[Belov \cite{be}]
Assume 
$ \displaystyle a_k \searrow 0. $
Then $\sum a_k\,\sin(kx) $ is a PS in $ [0,\pi ] $ iff

\begin{equation}  \sum_{k=1}^{n} (-1)^{k-1}\,ka_k \geq 0, \qquad  \mbox{for all } n\geq 2.  \Label{bel}  \end{equation}
\end{THEOREM}
}

\begin{REM} \em
For the cosine analog, condition (\ref{bel})
is sufficient but not necessary.
\end{REM}

Belov's Theorem leaves no more room for improvement, 
unless the $ \searrow $ assumption on $ a_k $ is lifted. 
In this less restrictive situation, (\ref{bel}) is no longer sufficient for PS
(it is still necessary).  It is not difficult to construct
examples of PS sine sums with non-monotone coefficients, as we will see in Section~2.
However, no useful general conditions applicable to non-monotone coefficients
are known until recently. 
In \cite{Kw2}, the following result was established.
{
\renewcommand{\theTHEOREM}{C}
\begin{THEOREM}
Vietoris' result remains valid when $ \searrow $ (still need (\ref{v}))
is relaxed to
$$  \frac{(2j-1)\sqrt{j+1}}{2j\sqrt{j}} \,\,a_{2j+1}\, \leq  \,a_{2j}, \qquad  j=1,2,\cdots .  $$
\end{THEOREM}
}
An example is given by the non-monotone sequence of coefficients:
\begin{eqnarray}
  && 1\,, \hspace*{6mm} \frac{1}{2} \,, \hspace*{7mm} \frac{1}{\sqrt{2}} \,, \hspace*{12mm} \frac{3}{4\sqrt{2}} \,, \hspace*{12mm} \frac{1}{\sqrt{3}} \,, \hspace*{12mm} \frac{5}{6\sqrt{3}} \,, \hspace*{10mm} \cdots  \Label{Cm} \\[1.2ex]
  = \hspace*{-4mm}  && 1\,, \hspace*{4mm} 0.5\,, \hspace*{4mm} 0.707\cdots \,, \hspace*{4mm} 0.530\cdots \,, \hspace*{4mm} 0.577\cdots \,, \hspace*{4mm} 0.481\cdots \,, \hspace*{6mm} \nonumber  \cdots 
\end{eqnarray}

An important tool used in the proof is the well-known Comparison Principle
(CP for short). It will continue to play an important role in this
paper. 

Since a non-zero scalar multiple of a PS is still a PS, we 
consider two sequences of coefficients equivalent if they only differ by a non-zero multiple.
We say that 

\par\vspace*{-4mm}\par
\begin{center}
$ \left\{ a_k\right\} \succeq\left\{ b_k\right\}  $  $ \Longleftrightarrow  $
(1)  $ a_k=0\Longrightarrow b_k=0 $, and
(2)  after skipping those $ a_k=0 $, $\,\,\displaystyle\frac{b_k}{a_k} \searrow\,0$.
\end{center}
\par\vspace*{-4mm}\par

\noindent
This defines a partial ordering among equivalent classes of sequences.
With this notation, the CP can be restated as follows.

\begin{LEMMA}\Label{cmp}
Let $ \sigma _k(x) $ be a sequence of functions defined on an
interval $ I $.
\begin{center}
     $ \sum a_k\sigma _k(x) $ PS in $ I $ and 
     $ \left\{ a_k\right\} \succeq\left\{ b_k\right\} \,\,\Longrightarrow \,\,\sum b_k\sigma _k(x) $ PS in $ I $.
\end{center}
\end{LEMMA}

\begin{REM} \em
Among all the sequences of coefficients satisfying Vietoris' conditions, there is a maximal
one, namely
\begin{equation}  \left\{ c_k\right\}  =  \left\{    1\,, \hspace*{6mm} \frac{1}{2} \,, \hspace*{6mm} \frac{1}{2} \,, \hspace*{6mm} \frac{3}{8} \,, \hspace*{6mm} \frac{3}{8} \,, \hspace*{6mm} \frac{5}{16} \,, \hspace*{6mm} \cdots  \right\}   \end{equation}
obtained by replacing the inequality sign in (\ref{v}) by equality and 
letting $ a_{2j}=a_{2j+1} $. The CP 
reduces the proof of the general Vietoris inequality to just
showing that the maximal sum $ \mb V=\sum\,c_k\,\sin(kx) $ is PS.

In the same sense, (\ref{Cm}) is the maximal sequence for Theorem~C. 
On the other hand, there is no maximal sequence for Belov's result.
\end{REM}

In this paper, we present two further improvements of Theorem~C,
In order to better illustrate some of the main ideas, we
first establish, in Section~\ref{sc3}, a slightly weaker NN criterion
that is associated with the sequence of coefficients
\begin{equation}  \left\{ \gamma _k\right\}  = \left\{  2, \, 1, \, \frac{4}{3} , \, 1, \, \frac{6}{5} , \, 1, \, \frac{8}{7} , \, 1 , \, \cdots \right\}  , \qquad  \gamma _k = \begin{cases} \mbox{\footnotesize  $\displaystyle \frac{k+1}{k} $} \quad  & k \mbox{ is odd} \\[1.5ex]  1 \quad  & k \mbox{ is even} \end{cases} \,\, .  \Label{ak1}  \end{equation}

\setcounter{THEOREM}{0}
\begin{LEMMA}
$ \mb\Psi =\sum a_k\sin(kx) $ is a PS in $ [0,\pi ] $ if (\ref{v}) holds and
\begin{equation}  a_{2j}\, \geq \, \frac{2j+1}{2j+2} \,\,a_{2j+1}, \qquad  \mbox{for } j=1,2,\cdots   .  \Label{kv}  \end{equation}
The maximal sum is given by $ \mb{\Phi }=\sum\,\gamma _k\,\sin(kx) $.
\end{LEMMA}

\par\vspace*{4mm}\par
\begin{REM} \em
Lemma 2 is already a significant improvement over Theorem~A and~C.
The coefficients $ a_k $ that satisfy the hypotheses of these Theorems
must decay faster than $ 1/\sqrt{k} $. The coefficients of $ \mb\Phi  $, on the other hand,
converge to $ 1 $.
\end{REM}

Lemma~2 can be sharpened in two different ways. Let 
$ \displaystyle \alpha \approx  0.78265213271\cdots $
be the second largest real root of the polynomial
\begin{equation}  54675\,{a}^{4}-2442195\,{a}^{3}+2182800\,{a}^{2}-115424\,a-96429 = 0.  \Label{alp}  \end{equation}

\par\vspace*{1mm}\par
\begin{THEOREM}
$ \mb\Psi =\sum a_k\sin(kx) $ is a PS in $ [0,\pi ] $ if $ \left\{ a_k\right\}  $ satisfies
(\ref{v}),
\begin{equation}  a_{2}\, \geq \, \frac{3\alpha }{4}  \, a_{3} , \quad  \mbox{and (\ref{kv}) \,\,for\,\, } j=2,3,\cdots.  \end{equation}
The maximal sum is $ \mb\Phi _1 $ with coefficients
$ \displaystyle \left\{  2\alpha  \,,\, \alpha   \,,\, \gamma _3 \,,\, \gamma _4 \,,\, \gamma _5 \,,\, \cdots \right\} . $

The value $ \alpha  $ is best possible; if it is replaced by any
smaller positive number, then $ \mb\Phi _1(5) $ is not NN.
\end{THEOREM}
\par\vspace*{1mm}\par

\begin{REM} \em
Note that even
though the coefficients of $ \mb\Phi _1 $ are not monotone, the subsequence of 
odd-order coefficients is decreasing, while the even-order coefficients are constant.
Contrast this with $ \mb\Phi _2 $ defined below.
Its subsequence of even-order coefficients is increasing.
\end{REM}

Let 
\par\vspace*{-9mm}\par
\begin{equation}  \left\{ \delta _k\right\}  = \left\{  3, \,  \frac{3}{2} , \, \frac{7}{3} , \, \frac{7}{4} , \, \frac{11}{5} , \, \frac{11}{6} , \, \cdots \right\} , \qquad          \delta _k = 2 - \frac{(-1)^k}{k} \, .  \end{equation}

\par\vspace*{4mm}\par
\begin{THEOREM}
$ \mb\Psi =\sum a_k\sin(kx) $ is a PS in $ [0,\pi ] $ if $ \left\{ a_k\right\}  $ satisfies (\ref{v}) and, 
\begin{equation}  a_{2j}\, \geq \, \frac{(2j+1)(4j-1)}{2j(4j+3)} \,\,a_{2j+1}, \qquad  \mbox{for } j=1,2,\cdots.  \Label{kv2}  \end{equation}
The maximal sum is $ \mb\Phi _2=\sum \delta _k\,\sin(kx) $.
\end{THEOREM}
\par\vspace*{1mm}\par

\begin{REM} \em
Theorems~1 and 2 are independent of each other 
as their extremal sums are not related to each other by $ \succeq $.
The same is true for Lemma 2 and Theorem~C.
On the other hand, each of Theorems~1 and~2 implies both Lemma 2 and Theorem~C. 
Yet, neither extends Belov's result.
It would be ideal
if Belov's result can be combined with Theorems~1 and~2 in a general unified way, 
but that remains a future goal for now. 
\end{REM}

By applying the reflection $ x\mapsto(\pi -x) $ to $ \mb\Phi  $ (or $ \mb\Phi _1 $
and $ \mb\Phi _2 $), we see that
its PS property is equivalent to that of
\begin{equation}  \mb{\Theta } =   \sum (-1)^{k+1}\,\gamma _k\sin(kx)  \Label{H}  \end{equation}
(and the corresponding $ \mb\Theta _i $, $ i=1,2 $). 

For any $ k\in(1,\infty ) $, define
\begin{equation}  \phi _k(x) = \sin((k-1)x) + \frac{k-1}{k} \, \sin(kx) \, ,  \Label{fk}  \end{equation}
\begin{equation}  \theta _k(x) = \sin((k-1)x) - \frac{k-1}{k} \, \sin(kx) \, .  \Label{hk}  \end{equation}

The partial sums $ \mb{\Phi }(n) $ and $ \mb{\Theta }(n) $ have the representations
\begin{equation}  \mb{\Phi }(n) = 2\,\phi _2(x) + \frac{4}{3} \,\phi _4(x) + \cdots + \frac{2\tilde n}{2\tilde n-1}  {\,\phi _{2\tilde n}(x)}  + \left[  \, \frac{(2\tilde n+2)\sin(nx)}{2\tilde n +1}  \,\right]  ,  \Label{Fj}  \end{equation}
\begin{equation}  \mb{\Theta }(n) = 2\,\theta _2(x) + \frac{4}{3} \,\theta _4(x) + \cdots + \frac{2\tilde n}{2\tilde n-1}  {\,\theta _{2\tilde n}(x)}  + \left[  \, \frac{(2\tilde n+2)\sin(nx)}{2\tilde n +1}  \,\right]  ,  \Label{Hj}  \end{equation}
where 
$ \tilde n  $
denotes the largest integer less than or equal to $ n/2 $, and the notation
$ \left[  \, \cdot \, \right]  $
means that the term is present only if $ n $ is an odd integer. 

\begin{REM} \em
An alternative way to see that Theorem~1 implies Lemma~2 is to note that
$$  \mb \Phi  = 2(1-\alpha )\phi _2 + \mb \Phi _1 .  $$
The first term on the righthand side is NN and the second term is a PS.
Likewise,
\begin{equation}  \mb \Phi  = \mb F + \mb \Phi _2 ,  \Label{Ph1}  \end{equation}
where $ \mb F $ is the Ferj\'er-Jackson-Gronwall PS,
shows that Theorem~2 implies Lemma~2.
\end{REM}

The following well-known identities will be used in subsequent proofs.

{\footnotesize 
\begin{eqnarray}
\sin(x)+\sin(3x)+\sin(5x)+\cdots+\sin((2n-1)x) \hspace*{1.4mm} &=& \frac{1-\cos(2nx)}{2\cos(x)} \,.  \Label{s3} \\
\sin(x)\,+\,\sin(2x)\,+\,\sin(3x)\,+\,\cdots\,+\,\sin((nx) \hspace*{5.8mm} &=& \frac{\cos(\frac{x}{2} )-\cos(\frac{(2n+1)x}{2} )}{2\sin(\frac{x}{2} )} \,.  \Label{s1} \\
\cos(x)+\cos(3x)+\cos(5x)+\cdots+\cos((2n-1)x) \hspace*{0.4mm} &=& \frac{\sin(2nx)}{2\sin(x)} \,.  \Label{c3} \\
\cos(x)-\cos(2x)+\cos(3x)-\cdots+(-1)^n\cos((nx) &=& \frac{1}{2} + (-1)^n \frac{\cos(\frac{(2n+1)x}{2} )}{2\cos(\frac{x}{2} )} \,.  \Label{c1}
\end{eqnarray}
}

The rest of the paper is organized as follows. In Section~\ref{sc2}, we give
some examples of PS sine sums with non-monotone coefficients that can be easily
constructed using known results. These examples should be contrasted 
with those covered by Theorems~1 and~2. The proofs of Lemma~2 and Theorems~1 and~2
are given in Sections 3, 4, and 5, respectively. Section~\ref{sc6} presents some
further examples and remarks.

\section{Trivial Examples of PS with Non-Monotone Coefficients \label{sc2}}

\begin{EX} \em
Assume $ b_k\searrow0 $. Then 
$ \mb{B} = \sum b_k\sin((2k-1)x) $
is a PS in $ [0,\pi ] $.
\end{EX}
Consider 
\begin{equation}  \mb{C}(n)=\sin(x)+\sin(3x)+ \cdots +\sin((2n-1)x) .  \end{equation}
From (\ref{s3}), we see that
\begin{equation}  2\cos(x) \, \mb{C}(n) = 1 - \cos(2nx) \geq  0 .  \end{equation}
Hence, $ \mb{C} $ is PS. It follows from the CP that
$ \mb{B} $ is also PS.

Even though the sequence $ \left\{ b_k\right\}  $ is decreasing, 
from the point of view of the
full sine sum, the coefficient sequence is actually $ \left\{ b_1,0,b_2,0,b_3,0,\cdots\right\}  $,
which is not monotone.

\par\vspace*{\baselineskip}\par

\begin{EX} \em Let $ \mb B $ and $ \mb C $ be as in Example 1 and $ \mb V $ be the Vietoris
sum as in Remark 3.
\begin{equation}  \mb{C}+\mb{V} = 2\sin(x)+\frac{1}{2} \,\sin(2x) + \frac{3}{2} \,\sin(3x) + \frac{3}{8} \,\sin(4x) + \cdots ,  \Label{SV}  \end{equation}
is a PS with non-monotone coefficients.
More generally, $ \beta \mb{B}+\mb{V} $ is a PS for any $ \beta >0 $.
\end{EX}

\par\vspace*{\baselineskip}\par

\begin{EX} \em
By applying the reflection  $ x\mapsto\pi -x $ to $ \mb{V} $, we see that 
\begin{equation}  \mb{V}_2 = \sum (-1)^{k+1} c_k\,\sin(kx)  \end{equation}
is a PS in $ [0,\pi ] $, so is
$ 2\mb{V} + \mb{V}_2  $
with coefficients
\begin{equation}  3\,, \hspace*{2mm} \,\frac{1}{2} \,, \hspace*{2mm} \,\frac{3}{2} \,, \hspace*{1mm} \,\frac{3}{8} \,, \hspace*{1mm} \,{\frac {9}{8}}\,, \hspace*{1mm} \, {\frac {5}{16}}\,, \hspace*{1mm} \, ...   \end{equation}
\end{EX}

\par\vspace*{\baselineskip}\par

\begin{EX} \em
It is easy to construct specific 
sine polynomials with a finite number 
of terms and non-monotone coefficients that are PS in $ [0,\pi ] $.
For example
\begin{equation}  2\sin(x) + \sin(2x) + \left( 1+\frac{\sqrt3}{2} \right) \,\sin(3x)   \Label{po1}  \end{equation}
and
\begin{equation}  3\sin(x) + \sin(2x) + \left( \frac{3}{2} + \sqrt2 \right) \,\sin(3x)   \Label{po2}  \end{equation}
are both PS in $ [0,\pi ] $ with non-monotone coefficients. We refer the readers to
\cite{Kw} for a discussion of how these and similar polynomials can be constructed.

It is also easy to prove that for any positive integer $ m $,
$$  \sin(x) + \frac{\sin(mx)}{m}   $$
is a PS in $ [0,;\pi ] $ with non-monotone coefficients.

If one insists on constructing examples with an infinite number of terms,
simply add an appropriate multiple of one of these to $ \mb{V} $.
\end{EX}

We consider all such examples trivial because they are easy corollaries
of Vietoris' result and other known examples.

\section{Proof of Lemma 2 \label{sc3}}

Lemma 2 is obviously true for $ n=1,2 $ and $ 3 $.
Hence, we assume $ n\geq 4 $ in the following.

\begin{LEMMA}\Label{jk}
For all $ k>1 $,
\begin{equation}  \theta _k(x)\geq 0 \quad  \mbox{ for } x \in \left[ 0,\frac{\sigma }{k} \right]  ,  \Label{jp}  \end{equation}
where $ \sigma \approx 4.493409458 $ is the first positive zero of the function
\begin{equation}  f(z) = \sin(z) - z\cos(z) .  \Label{f}  \end{equation}
\end{LEMMA}  \begin{PROOF}
Let 
$ \displaystyle \mu =1-\frac{1}{k} \in \left( 0 ,1\right)  $
and $ y=kx $. Then, from the definition (\ref{fk}),
\begin{equation}  \frac{\theta _k(x)}{\mu } = \frac{\sin(\mu y)}{\mu } - \sin(y) .  \end{equation}
\begin{equation}  \frac{\partial }{\partial \mu } \left(  \frac{\theta _k(y)}{\mu } \right)  = - \frac{\sin(\mu y)-\mu y\cos(\mu y)}{\mu ^2} = - \frac{f(\mu y)}{\mu ^2} \,.  \Label{jpa}  \end{equation}
For $ x\in[0,\sigma /k] $, $ \mu y\in[0,\sigma ] $. Since $ f(z) $ is positive in $ (0,\sigma ) $, the 
righthand side of (\ref{jpa}) is negative, implying that $ \theta _k(y)/\mu  $
is a decreasing function of $ \mu  $. Hence,
\begin{equation}  \frac{\theta _k(x)}{\mu } \geq  \lim_{k\rightarrow \infty } \frac{\theta _k(x)}{\mu } = 0 .  \end{equation}
\end{PROOF}

\begin{LEMMA}\Label{H1}
For any integer $ n>0 $,
\begin{equation}  \mb{\Phi }(n)\geq 0 \quad  \mbox{in }  \left[ 0, \frac{\pi }{n} \right]  \cup \left[  \pi -\frac{\pi }{n} \, , \pi \right]  .  \end{equation}
\end{LEMMA}  \begin{PROOF}
In $ [0,\pi /n] $, every term in $ \mb\Phi (n) $ is NN 
and so is their sum.

The assertion $ \mb\Phi (n)\geq 0 $ in $ [\pi -\pi /n,\pi ] $ is equivalent to $ \mb\Theta  $
being NN in $ [0,\pi /n] $. We make use of the representation (\ref{Hj}) of $ \mb\Theta (n) $.
If $ n $ is even, $ \mb\Theta  $ is a sum of positive multiples of $ \theta _{2j}(x) $,
for $ j=1,\cdots,\tilde n $.
By Lemma~\ref{jk}, each of these is NN in $ [0,\sigma /2\tilde n]\supset[0,\pi /n] $.
Hence, their sums is
NN in $ [0,\pi /n] $. If $ n $ is odd, 
there is an extra term $ \sin(nx) $ which is
also NN in $ [0,\pi /n] $ and the conclusion still holds.
\end{PROOF}

In view of Lemma~\ref{H1}, to complete the proof of Lemma~2,
it remains to show that $ \mb\Phi (n) $ is
NN in $ I_n=[\pi /n,\pi -\pi /n] $ for all $ n $.

Let $ m=n $ if $ n $ is odd, and $ n-1 $ otherwise. It is the 
largest odd integer $ {}\leq n $.
Then $ \mb\Phi (n)=\mb{S}_1(n)+\mb{T}(m) $, where
\begin{equation}  \mb S(n) = \sin(x)+\sin(2x)+\cdots+\sin(nx)   \Label{S1}  \end{equation}
and
\begin{equation}  \mb{T}(m) = \sin(x)+\frac{\sin(3x)}{3} + \cdots + \frac{\sin(m x)}{m}.  \Label{s2}  \end{equation}
Identity (\ref{s1}) gives a lower bound for $ \mb S(n) $.
\begin{eqnarray}
  \mb S(n) &\geq & \frac{\cos(x/2)-1}{2\sin(x/2)} \nonumber \\[1.2ex]
	  &=& - \frac{\tan(x/4)}{2} \Label{Stan} \\[1.2ex]
	  &\geq & - \frac{1}{2} \,.
\end{eqnarray}
% It is easy to show that $ \phi _2(x)\leq \phi _2(2\pi /3)=3\sqrt3/4 $.
% Hence,
The proof of Lemma 2 is thus complete if we can show that
\begin{equation}  \mb{T}(m) \geq  \frac{1}{2} \, , \qquad  x \in I_n, \quad  n\geq 4 .  \Label{S23}  \end{equation}

When $ n $ is even, $ n $ and $ n-1 $ use the same
$ \mb{T}(m) $, but $ I_{n-1}\subset I_n $.
Hence, if (\ref{S23}) can be proved for $ n $, then it will
also hold for $ n-1 $. In other words, we only have to establish (\ref{S23})
for even $ n $, in which case $ m=n-1 $.
Note that $ \mb{T}(m) $ is an even function about $ x=\pi /2 $. Thus, 
it suffices to show (\ref{S23}) for odd $ m $ and $ x\in J_n=[\pi /n,\pi /2] $.

An alternative representation for $ \mb{T}(m) $ can be given using (\ref{c3}).
\begin{eqnarray}
  \mb{T}(m) = f_n(x)  &:=& \int_{0}^{x} \left(  \cos(s)+\cos(3s) + \cdots + \cos((n-1) s) \right)  \,ds \nonumber \\[1.2ex]
      &=& \int_{0}^{x} \frac{\sin(ns)}{2\sin(s)} \,ds  \, .    \Label{sx}
\end{eqnarray}
For convenience, we revert back to using $ n=m-1 $ instead of $ m $.
Besides being easier to estimate, another advantage of the alternative 
representation is that the definition of $ f_n(x) $ can be extended to all
real $ n\in(0,\infty ) $. Even though we only need
(\ref{S23}) for even integers $ n $, we are going to prove the stronger inequality
\begin{equation}  f_n(x) \geq  \frac{1}{2} \, , \qquad  x \in J_n, \quad  n\geq 4 .  \Label{S24}  \end{equation}

\begin{center}
\FG{59}{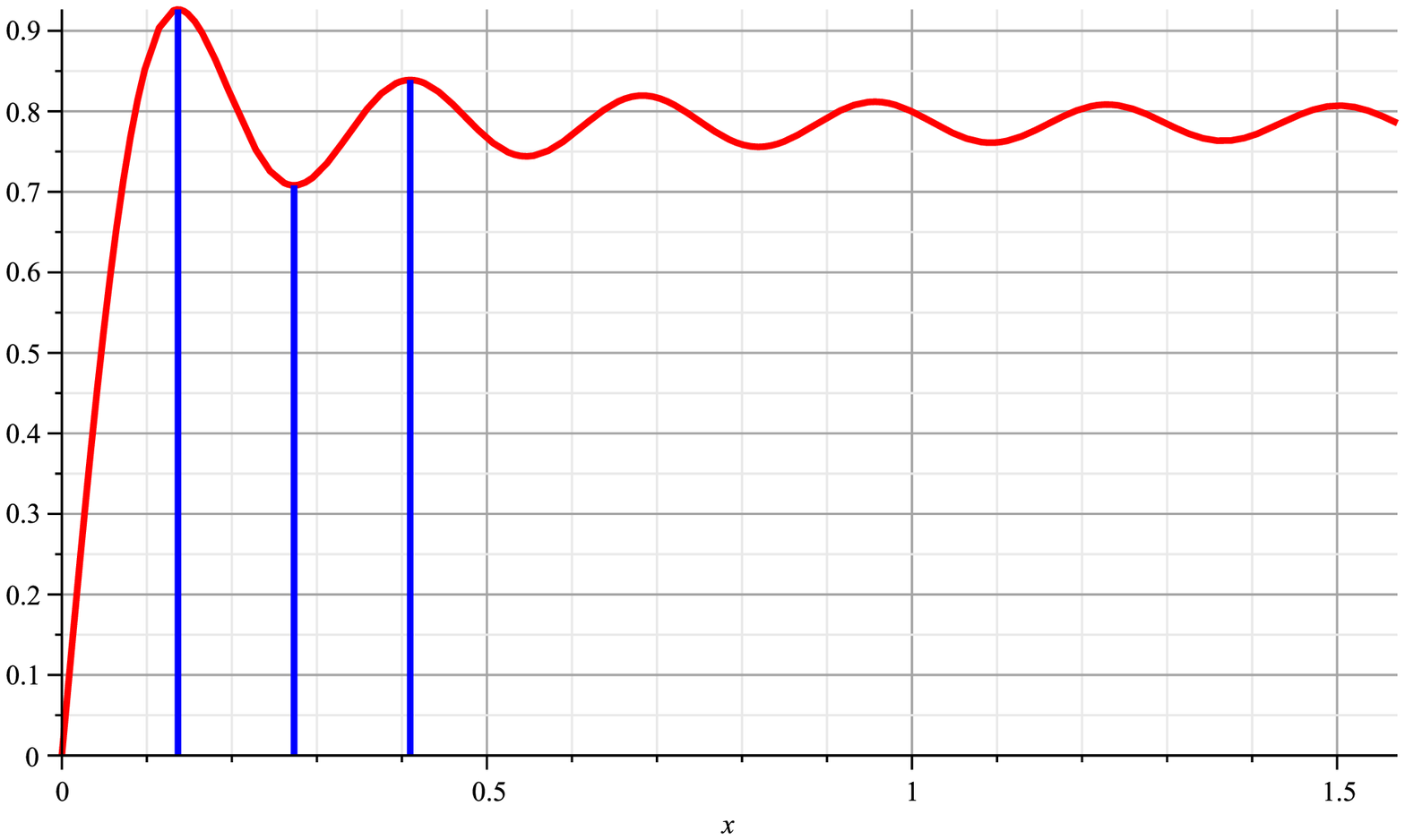} \\
\par\vspace*{-7.8mm}\par
{\footnotesize  \hspace*{4.0mm} $ x_1 $\hspace*{3.6mm} $ x_2 $\hspace*{3.6mm} $ x_3 $\hspace*{58mm} }
\par\vspace*{1mm}\par
Figure 1. Graph of $ f_{23}(x) $.
\end{center}

\par\vspace*{3mm}\par

The graph of one of these functions, $ f_{23}(x) $, is depicted in Figure 1.

Since $ f_n'(x)=\sin(nx)/\sin(x) $, the critical
points of $ f_n(x) $ in $ J_n $ are  $ \pi /n,\,2\pi /n,\,3\pi /n,\,\cdots $.
The first is the left endpoint of $ J_n $ and is a local maximum,
so are all other odd-order points.
The even-order points  $ x_2=2\pi /n,\,x_4=4\pi /n,\,\cdots $ are
local minima. 
A lower bound for $ f_n(x) $ in $ J_n $ is, therefore,
\begin{equation}  \min_{x\in J_n}  f_n(x) = \min \left\{ f_n(x_2), f_n(x_4), \cdots \right\} .  \Label{s4}  \end{equation}

Integration by parts gives
\begin{eqnarray}
   f_n(x_{2k+2}) - f_n(x_{2k}) &=& \int_{x_{2k}}^{x_{2k+2}} \frac{\sin(ns)}{2\sin(s)} \,ds  \nonumber \\[1.2ex]
            &=&  \int_{x_{2k}}^{x_{2k+2}} \frac{(1-\cos(ns))\cos(s)}{2n\sin^2(s)} \,ds  \nonumber \\[1.2ex]
	    &>& 0 . \Label{s6}
\end{eqnarray}

Hence, $ f_n(x_2)<f_n(x_4)<f_n(x_6)<\cdots $ and it follows from (\ref{s4}) that
\begin{equation}  f_n(x) \geq  f_n(x_2) .  \Label{s5}  \end{equation}
Now (\ref{S24}) follows from the next Lemma and the proof of Lemma~2 is complete.

\begin{LEMMA}\Label{S2}
The sequence $ f_n(x_2) $, $ n=4,5,\cdots $ is increasing.
\begin{equation}  \frac{2}{3} = f_4(x_2) < f_5(x_2) < ... < f_n(x_2) < f_{n+1}(x_2) < ...  \end{equation}
\end{LEMMA}  \begin{PROOF}
That $ f_4(x_2)=2/3 $ can be verified directly. In fact, each $ f_n(x_2) $ can
be computed exactly using Maple.

The change of variable, $ s=t/n $ gives
\begin{equation}  f_n(x_2) = \int_{0}^{2\pi /n} \frac{\sin(ns)}{2\sin(s)} \,ds = \int_{0}^{2\pi } \frac{\sin(t)}{2n\sin(t/n)} \,dt  = \int_{0}^{\pi } k_n(t)\sin(t)\,ds,  \end{equation}
where
\begin{equation}  k_n(t) = \frac{1}{2n\sin(t/n)} \,.  \end{equation}
Thus,
\begin{equation}  f_{n+1}(x_2)-f_n(x_2) = \int_{0}^{2\pi } \big( k_{n+1}(t) - k_n(t) \big) \sin(t) \,dt .  \end{equation}
In the next Lemma, we show that 
\begin{equation}  h_n(t)=k_{n}(t)-k_{n+1}(t)   \Label{ht}  \end{equation}
is a positive increasing function of $ t\in[0,2\pi ] $. Anticipating this fact, we see that
\begin{eqnarray}
    f_{n+1}(x_2)-f_n(x_2) &=& \int_{\pi }^{2\pi } |\sin(t)|h_n(t) \,dt -  \int_{0}^{\pi } \sin(t)h_n(t) \,dt \nonumber \\
     &>& h_n(\pi ) \int_{\pi }^{2\pi } |\sin(t)| \,dt - h_n(\pi ) \int_{0}^{\pi } \sin(t) \,dt  \nonumber \\
     &=& 0.
\end{eqnarray}
as desired.
\end{PROOF}

\begin{LEMMA}
For $ n\geq 4 $, $ h_n(t) $ is a positive increasing function of $ t $ in $ [0,2\pi ] $.
\end{LEMMA}  \begin{PROOF}
The NN of $ h_n(t) $ follows from the fact that,
for fixed t, $ k_n(t) $ is a decreasing
function of $ n $, which is equivalent to the fact that $ n\sin(t/n) $ is 
an increasing function of $ n $.

The increasing property of $ h_n(t) $ is true if we can prove that
\begin{equation}  \frac{\partial ^2}{\partial n\partial t} \, k_n(t) \leq  0.   \Label{hnt}  \end{equation}
Direct computation gives the numerator of 
$ -\,\frac{\partial ^2}{\partial n\partial t} \, k_n(t)  $
as the function
\begin{equation}  \xi (t) = 2\cos^2(t/n) +t\,\sin^2(t/n) - 2n \cos(t/n) \sin(t/n) .  \end{equation}
For convenience, we have suppressed the dependence of $ \xi (t) $ on $ n $.
Now it suffices to show that $ \xi (t)\geq 0 $ for $ t\in[0,2\pi ] $. Since
$ \xi (0)=0 $, if we can show that $ \xi '(t)\geq 0 $,
the proof is complete.
\begin{equation}  \xi '(t) = \sin\left(  \frac{t}{n} \right)  \left[  3 \sin \left(  \frac{t}{n} \right)  - \frac{2t}{n} \,\cos\left( \frac{t}{n} \right)   \right]   = \sin\left(  \frac{t}{n} \right) \xi _2(t).  \end{equation}
It now suffices to show that $ \xi _2(t) $
is NN. The desired conclusion follows
from the facts $ \xi _2(0)=0 $, and 
\begin{equation}  \xi _2'(t) = \frac{1}{n} \, \cos\left( \frac{t}{n} \right)  + \frac{2t}{n^2} \sin\left( \frac{t}{n} \right)  \geq  0 \,.  \end{equation}
\end{PROOF}

\section{Proof of Theorem 1 \label{sc4}}

We first take care of $ n>20 $.
The partial sums $ \mb\Phi _1(n) $ can be represented as
\begin{eqnarray}
            \mb{\Phi }_1(n) &=& \mb\Phi (n) - \lambda  \phi _2(x)   \nonumber \\
	                 &=& \mb{S}(n) +  f_n(x)  - \lambda  \phi _2(x), 
\end{eqnarray}
where $ \lambda =2-2\alpha \approx 0.434695735 $.  In view of (\ref{Stan}) and Lemma~\ref{S2}, we get,
for all $ n>20 $, $ x\in[0,\pi ] $,
\begin{equation}  \mb{\Phi }_1(n) \geq  F(x) := - \,\frac{\tan(x/4)}{2} + f_{20}(x_2) -\frac{4347}{10000} \, \phi _2(x).  \Label{F20}  \end{equation}
Maple gives
\begin{equation}  f_{20}(x_2) =  \frac{2}{15} +{\frac {1580}{4641}}\,\cos \left( \frac{\pi }{5}  \right) +{\frac {1820 }{1881}}\,\cos \left(\frac{2\pi }{5}  \right) > \frac{73542}{103909} \, .  \Label{F20a}  \end{equation}
It follows from (\ref{F20}) and (\ref{F20a}) that
\begin{equation}  \mb{\Phi }_1(n) \geq  F_1(x) := \frac{73542}{103909}  - \,\frac{\tan(x/4)}{2}  -\frac{4347}{10000} \left(  \sin(x) + \frac{\sin(2x)}{2} \right) .  \Label{F2a}  \end{equation}
Let $ T=\tan(x/4) $.
Since $ x\in[0,\pi ] $, we have $ T\in[0,1] $. Then
\begin{eqnarray*}
   F_1(x) &=& {\frac {73542}{103909}} - \frac{T}{2} -{\frac {4347}{1250}}\,{\frac {T \left( 1-{T}^{2} \right) ^{3}}{ \left( 1+{T}^{2} \right) ^{4}}} \\[1.5ex]
          &=& \frac{P(T)}{(1+T^2)^4} \,.  
\end{eqnarray*}
where
{\small 
\begin{eqnarray*}
   P(T) \!\! &=& \!\! -45963750\,{T}^{9}+91927500\,{T}^{8}+267837423\,{T}^{7}+367710000\,{T} ^{6}-1630859769\,{T}^{5} \\
   &&  +551565000\,{T}^{4}+1171222269\,{T}^{3}+ 367710000\,{T}^{2}-497656173\,T+91927500 .
\end{eqnarray*}
}

The classical Sturm Theorem, provides a way to find the 
number of real roots of an algebraic polynomial with real coefficients 
within any given subinterval of the real line. It can be used (see \cite{Kw} and the
discussion below) to show that $ P(T)>0 $ in $ [0,1] $.
With this fact, we conclude that $ \mb\Phi _1(n)>0 $ for $ x\in[0,\pi ],\,n>20 $.

For $ n\leq 20 $, the above argument does not work because when $ f_{20}(x_2) $ is
replaced by any $ f_n(x_2) $ with $ n<20 $, the resulting $ F(x) $ is no longer
NN in $ [0,\pi ] $.
Our verification of Theorem~1 for $ n\leq 20 $, relies on a brute-force
technique based on the Sturm Theorem.
The method is explained in great details in \cite{Kw}.
See also \cite{AK} which discusses its use in 
the study of Rogosinski-Szeg\"o-type inequalities \cite{AK2}.

In a nutshell, given any specific sine polynomial, we can expand it into a product
of $ \sin(x) $ and an algebraic polynomial $ p(Y) $ of the variable $ Y=\cos(x)\in[-1,1] $.
The Sturm Theorem can then be invoked to check if $ p(Y) $ is NN or not. 

This procedure works with one polynomial at a time. It is, therefore, not
adequate to prove general results like Theorem~1, which involves an infinite number 
of polynomials. Nevertheless, we can comfortably use this technique to deal with the 
first 20 of such polynomials. 

The procedure we implemented in
Maple, however, has one limitation. It works only
when the coefficients of the sine polynomial are given rational numbers. For this reason,
it cannot be directly applied to the sine polynomials of Theorem~1 
because they involve the irrational number $ \alpha  $.
The procedure is modified as follows. For $ n\leq 20 $, except $ n=5 $, we replace
$ \alpha  $ by the slightly smaller rational number $ \underline{\alpha }=171/100<\alpha  $.
The corresponding partial sums $ \underline{\mb\Phi }_1(n) $ is shown to be
NN using the Maple procedure.  It then follows from the CP that $ \mb\Phi _1(n) $ is also NN.

With $ \mb\Phi _1(5) $, the above approach encounters
a different problem. No matter what $ \underline\alpha <\alpha  $ is chosen, 
$ \underline{\mb\Phi }_1 $ is not NN.
In fact, $ \alpha  $ has been chosen to be critical in some sense, namely,
$$  \alpha  = \inf \left\{ a \,\left |\, p_a(Y) \geq 0 \mbox{ in } [0,\pi ] \right .\right\}  .  $$
Here $ p_a(Y) $ is the algebraic polynomial
\begin{equation}  p_a(Y) = 144\,{Y}^{4}+60\,{Y}^{3}-68\,{Y}^{2}+ \left( 15\,a-30 \right) Y+(15\,a-1) .  \end{equation}
associated with the sine polynomial
\begin{equation}  2a  \,\sin(x) + a  \,\sin(2x) + \sum_{k=3}^{5} \gamma _k \,\sin(kx)\geq 0  .  \end{equation}
For large $ a $, for example $ a=2 $, $ p_a(Y) $ is NN
in $ [-1,1] $; its graph lies above and away from the $ Y $-axis. 
On the other hand, when $ a=0 $, the graph crossed the $ Y $-axis.
As $ a $ increases from 0, the graph of $ p_a(Y) $ rises monotonically.
By continuity, there is a value of $ a=\alpha  $ when the graph is about to
leave the $ Y $-axis; 
it is tangent to the $ Y $-axis at one or more points.
To determine $ \alpha  $, note that each point
of tangency corresponds to
a double root of $ p_\alpha (Y)=0 $. A necessary condition for having a double root is
the vanishing of the discriminant.
With the help of Maple, the discriminant, after deleting a numerical
factor, is found to be (\ref{alp}). Numerical computation yields four real 
roots of (\ref{alp}):
$ \displaystyle -0.17, \, 0.30, \, 0.78, $
and $ 43.76 $. Hence, $ \alpha  $ is the second largest root. 

\section{Proof of Theorem 2 \label{c5}}

As in the proof of Theorem 1, we can use the Sturm procedure to confirm 
Theorem~2 for small $ n $, more specifically, we have done that
for $ n\leq 20 $. Hence, we assume $ n>20 $ in the
rest of this section.

The partial sums of $ \mb\Phi _2 $ have the representation
\begin{equation}  \mb \Phi _2(n) = 2\mb S(n) + \mb U(n) ,  \Label{F2}  \end{equation}
where $ \mb{S} $ is given by (\ref{S1}) and
\begin{eqnarray}
     \mb U(n) &=& \sin(x) - \frac{\sin(2x)}{2} + \cdots - \frac{(-1)^n\sin(nx)}{n} \nonumber \\
            &=& \int_{0}^{x} \big(\cos(s) - \cos(2s) + \cdots - (-1)^n\cos(ns) \big) \,ds  \nonumber \\
	    &=&  \frac{x}{2} +(-1)^n \int_{0}^{x} \frac{\cos\big(\frac{(2n+1)s}2\big)}{2\cos\big(\frac{s}2\big)} \,ds .    \Label{Un}
\end{eqnarray}
We have used (\ref{c1}) to derive the last equality. By Lemma~\ref{H1}, we only have 
to show that $ \mb\Phi _2(n)\geq 0 $ in $ I_n=[\pi /n,\pi -\pi /n] $. Using (\ref{Stan}),
(\ref{F2}) and (\ref{Un}),
we see that
\begin{equation}  \mb\Phi _2(n) \geq  -\tan\left( \frac{x}{4} \right)  + \frac{x}{2} - h_n(x) ,  \Label{F2aa}  \end{equation}
where
\begin{equation}  h_n(x) = (-1)^{n+1} \int_{0}^{x} \frac{\cos\big(\frac{(2n+1)s}2\big)}{2\cos\big(\frac{s}2\big)} \,ds .  \end{equation}
Since $ \tan(x/4)\leq 0.32x $ for $ x\in[0,\pi ] $, (\ref{F2aa}) leads to
\begin{equation}  \mb\Phi _2(n) \geq  0.18 \,x  - h_n(x) .  \Label{F2b}  \end{equation}
Hence, Theorem~2 is proved if we can show that
\begin{equation}  h_n(x)  \leq  0.18\,x , \quad  x\in I_n, \quad  n\geq 21.  \Label{hn}  \end{equation}

With change of variables, $ s\mapsto 2t $, $ x\mapsto 2y $
and $ 2n+1\mapsto \hat m $, (\ref{hn}) 
becomes
\begin{equation}  \hspace*{32mm} g_{\hat{m}}(y)  \leq  0.18\,y ,  \qquad   y\in I_{\hat m}, \quad  \hat m=43,45,47,\cdots,  \Label{gm1}  \end{equation}
where
\begin{equation}  g_{\hat{m}}(y) = (-1)^{(\hat{m}+1)/2} \int_{0}^{y} \frac{\cos(\hat{m}t)}{2\cos(t)} \,dt   \Label{gm3}  \end{equation}
and $ I_{\hat m}=[\pi /(\hat m-1),\pi /2-\pi /(\hat m-1)] $.
In fact, we claim that (\ref{gm1}) holds in the bigger interval 
$ J_m=[\pi /\hat{m},\pi /2] $.

\begin{center}
\FG{59}{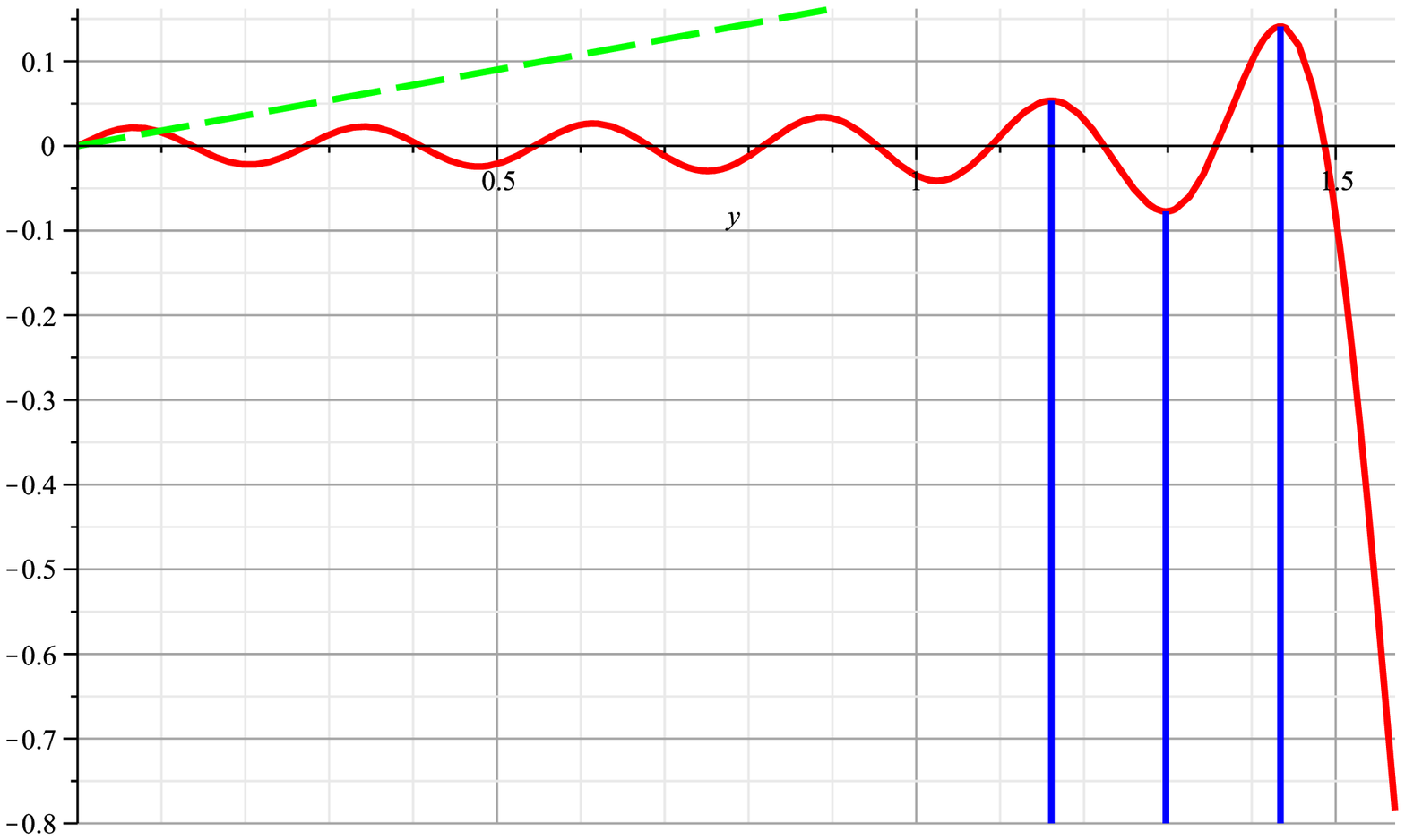} \\
\par\vspace*{-4.0mm}\par
{\footnotesize  \hspace*{65mm} $ y_3 $\hspace*{4mm} $ y_{2} $\hspace*{4mm} $ y_{1} $}
\par\vspace*{1mm}\par
\par\vspace*{1mm}\par
Figure 2. Graphs of $ u=g_{23}(y) $ and $ u=0.18y $.
\end{center}

\par\vspace*{3mm}\par

The wavy curve in Figure 2 depicts the graph of $ g_{23}(y) $ and the dashed
line is the graph of $ 0.18y $. It is clear from the figure that,
in this case, (\ref{gm1}) fails for small positive $ y $.
When $ \hat{m}=1(\mbox{mod}4) $, however, $ g_{\hat{m}}(y) $ is negative
for $ y\in[0,\pi /\hat{m}] $ and it can be shown that (\ref{gm1}) holds in the
whole interval $ [0,\pi /2] $.

The shape of $ g_{\hat{m}}(y) $ is strikingly similar to
that of $ f_m(x) $ in Figure 1.
Indeed, by using the reflection mapping $ y=\pi /2-x $, one can 
show that $ g_{\hat{m}}(y)=f_{\hat{m}}(\pi /2-y)-f_{\hat{m}}(\pi /2) $.
With this observation, we can
deduce many of the properties of $ g_{\hat{m}}(y) $ from those
of $ f_n(x) $.

For instance, the critical points of $ g_{\hat{m}}(y) $ are given by the sequence
$$  y_{(\hat{m}-1)/2}= \frac{\pi }{2\hat{m}} \quad  < \quad  y_{(\hat{m}-3)/2}=\frac{3\pi }{2\hat{m}} \quad  < \quad  \cdots \quad  < \quad  y_1= \frac{(\hat{m}-2)\pi }{2\hat{m}}  \,.  $$
Note that we have numbered the critical points $ y_k $ from right to left.
The first one, $ y_1 $, is always a local maximum and 
then they alternate as local minimum and maximum.
The last one, $ y_{(\hat{m}-1)/2} $ is is a minimum or
maximum depending on whether $ (\hat{m}+1)/2 $ is odd or even.
The sequence of local maximum (minimum) values $ g_{\hat{m}}(y_i) $ is decreasing 
(increasing) as $ i $ increases.
The global maximum of $ g_{\hat{m}}(y) $ is attained at $ y_1 $.

\begin{LEMMA}\Label{gm}
For all odd integers $ \hat{m}\geq 43 $,
\begin{equation}  g_{\hat{m}}(y) \leq   0.22, \quad  y\in[0,\pi /2].  \Label{gmy}  \end{equation}
\end{LEMMA}  \begin{PROOF}
Let us estimate 
\begin{eqnarray*}
       g_{\hat{m}}(y_1)-g_{\hat{m}}(y_2) &=& \int_{y_2}^{y_1} \frac{|\cos(\hat{m}t)|}{2\cos(t)} \,dt \, \\[1.5ex]
        &=&  \int_{\pi }^{2\pi } \frac{|\sin(s)|}{2\hat{m}\,\sin(s/{\hat{m}})} \,ds \,.
\end{eqnarray*}
For fixed $ s\in[\pi ,2\pi ] $, $ \hat{m}\sin(s/\hat{m}) $ is an increasing function 
of $ \hat{m}\geq 43 $. As a result, $ g_{\hat{m}}(y_1)-g_{\hat{m}}(y_2) $ is a decreasing function of $ \hat{m} $. In particular,
\begin{equation}  g_{\hat{m}}(y_1)-g_{\hat{m}}(y_2) \leq  g_{43}(y_1)-g_{43}(y_2) = 0.21731814075\cdots .  \end{equation}
Here we have abused the notation: $ y_1 $ and $ y_2 $ on the lefthand side of the
inequality are different from those on the other side.
Since $ g_{\hat{m}}(y_2)<0 $, the desired conclusion follows.
\end{PROOF}

Obviously, Lemma~\ref{gm} implies that (\ref{gm1}) holds on $ [11/9,\pi /2] $.
It remains to show (\ref{gm1}) on $ [\pi /\hat{m},11/9] $. Our next Lemma shows that 
in this subinterval, (\ref{gmy}) can be greatly improved.

\begin{LEMMA}\Label{gm2}
For all odd integers $ \hat{m}\geq 43 $,
\begin{equation}  g_{\hat{m}}(y) \leq   0.06, \quad  y\in[0,11/9].  \Label{gmy2}  \end{equation}
\end{LEMMA}  \begin{PROOF}
For $ \hat{m}=43 $, the first (counting from $ y_1) $ local maximum that falls within
the subinterval $ [0,11/9] $ is $ y_5 $, and we compute
\begin{equation}  g_{43}(y_5) - g_{43}(y_6) = 0.059552923006\cdots .  \end{equation}
Using the same arguments as in the proof of Lemma~\ref{gm}, we see that
$ g_{\hat{m}}(y_5)-g_{\hat{m}}(y_6) $ is a decreasing function of $ \hat{m} $. Hence,
\begin{equation}  g_{{\hat{m}}}(y_5) - g_{{\hat{m}}}(y_6) < 0.059552923006\cdots   \end{equation}
and the desired conclusion follows.
\end{PROOF}

Lemma~\ref{gm2} implies that (\ref{gm1}) holds on $ [1/3,11/9] $.
It remains to show (\ref{gm1}) on $ [\pi /\hat{m},1/3] $. We use a different method to
estimate $ g_{\hat{m}}(y) $ in this interval. For $ t\in[0,1/3] $,
\begin{equation}  1 \leq  \frac{1}{\cos(t)} \leq  \frac{1}{\cos(1/3)} < 1.06 .  \end{equation}
It follows that
\begin{equation}  \frac{\cos(\hat{m}t)}{2\cos(t)} \leq   \frac{1}{2} \,\cos(\hat{m}t) + 0.03  \Label{gm4}  \end{equation}
and
\begin{equation}  - \, \frac{\cos(\hat{m}t)}{\cos(t)} \leq  -  \frac{1}{2} \,\cos(\hat{m}t) + 0.03 .   \Label{gm5}  \end{equation}

We consider two cases. When $ (\hat{m}+1)/2 $ is even, then from (\ref{gm3}) and (\ref{gm4}), 
we obtain 
\begin{equation}  g_{\hat{m}}(y) \leq  \frac{\sin(\hat{m}y)}{2\hat{m}} + 0.03y .  \end{equation}
It is not difficult to see that this implies (\ref{gm1}) in $ [\pi /\hat{m},1/3] $.

In the complementary case, when $ (\hat{m}+1)/2 $ is odd, we use (\ref{gm3}) and (\ref{gm5})
to obtain
\begin{equation}  g_{\hat{m}}(y) \leq  - \,\frac{\sin(my)}{2\hat{m}} + 0.03y .  \end{equation}
This implies (\ref{gm1}) in $ [0,1/3]\supset[\pi /\hat{m},1/3] $, and completes the proof of
Theorem~2.

\section{Further Examples and Remarks \label{sc6}}

\begin{EX} \em
Theorem~1 can be applied to show that the sum
$$  \frac{\phi _2(x)}{\sqrt{2}} + \frac{\phi _4(x)}{\sqrt{3}} + \frac{\phi _6(x)}{\sqrt{4}} \cdots + \frac{\phi _{2\tilde{n}}(x)}{\sqrt{\tilde n +1}}  + \left[  \, \frac{\sin(nx)}{\sqrt{\tilde n+2}}  \,\right]    $$
is a PS. It is not covered by Theorem~2. More generally, Theorem~1 implies that
$$  \frac{\phi _2(x)}{\sqrt{\beta +1}} + \frac{\phi _4(x)}{\sqrt{\beta +2}} + \frac{\phi _6(x)}{\sqrt{\beta +3}} \cdots + \frac{\phi _{2\tilde{n}}(x)}{\sqrt{\beta +\tilde n}}  + \left[  \, \frac{\sin(nx)}{\sqrt{\beta +\tilde n+1}}  \,\right]    $$
is PS for $ \beta \geq \displaystyle\frac{8-9\alpha ^2}{9\alpha ^2-4} \approx 1.64393 $. Numerical experiments
suggest that the sum is PS for $ \beta >1.76923 $.
\end{EX}

\par\vspace*{\baselineskip}\par

\begin{EX} \em
Theorem~1 implies that
$$  \phi _2(x) + \frac{\phi _4(x)}{2^\gamma } + \cdots + \frac{\phi _{2\tilde{n}}(x)}{{\tilde n}^\gamma }  + \left[  \, \frac{\sin(nx)}{(\tilde n + 1)^\gamma }  \,\right]    $$
is PS for $ \gamma \geq 0.26 $. Theorem~2 performs worse in this case, giving only 
$ \gamma \geq 0.36258 $.
Numerical experiments suggest that the sum
may be a PS for $ 0.24\leq \gamma <0.26 $, but not for $ \gamma =0.23 $.
In the latter case, 
all partial sums except the sixth are NN in $ [0,\pi ] $.

These two examples indicate that Theorem~1 and 2 are not best possible.
\end{EX}

\par\vspace*{2mm}\par

\begin{REM} \em
In Theorems A, C, 1  and 2, the extremal sums are characterized by their
respective subsequences of odd-order coefficients, namely
$$  \left\{ c_{2j-1}\right\}  = \left\{ 1,\,\frac{1}{2} ,\,\frac{3}{8} ,\,\frac{5}{16} ,\,\cdots  \right\} ,  $$
$$  \hspace*{17mm} \left\{  1,\, \frac{1}{\sqrt2} ,\, \frac{1}{\sqrt3} ,\, \cdots  \right\} ,  $$
$$  \left\{  2\alpha , \gamma _{2j+1} \right\}  = \left\{  2\alpha , \,\frac{4}{3} ,\, \frac{6}{5} ,\, \frac{8}{7} ,\,\cdots  \right\} , \hspace*{7mm}  $$
and
$$  \hspace*{1mm} \left\{ \delta _{2j-1}\right\}  = \left\{  3, \,\frac{7}{3} ,\, \frac{11}{5} ,\, \frac{15}{7} ,\,\cdots  \right\} .  $$
The relative strength of the various results can be determined by
comparing these sequences according to the CP. For instance, 
sequence 1 $ \succeq $ sequence 2, while each of sequences 3 and 4 is $ \succeq $
sequences 1 and 2. To look for an improvement of Theorems~1 and 2, one searches
find a sequence $ \succeq $ sequence 3 or 4 that yields a PS. Note that 
$ \left\{ 1,1,\cdots \right\}  $ $ \succeq $ sequence 3 and 4, but its associated sine sum is not a PS.
In other words, $ \left\{ 1,\cdots \right\}  $ is a strict upper bound of all 
possible improvements of
Vietoris' sine result.
\end{REM}

\par\vspace*{2mm}\par

\begin{REM} \em
Theorem~1 relaxes the first condition, in (\ref{v}), of the 
Vietoris result. It is natural to ask
whether the second condition in (\ref{v}) can also be relaxed by replacing some
of the factors
$ \rho _j=\frac{2j-1}{2j} $ with larger constants. The following observation concerning
Belov's necessary condition (\ref{bel}) leads to the answer no.
\end{REM}

\begin{LEMMA}
{\rm (i)} A necessary condition for any sine polynomial 
$ \sum_{k=1}^{n} a_k\,\sin(a_kx) $ to be NN in some
neighborhood $ [\pi -\epsilon ,\pi ] $, $ 0<\epsilon <\pi  $ is 
\begin{equation}  \sum_{k=1}^{n} (-1)^{k-1}\,ka_k \geq 0.  \Label{bel1}  \end{equation}

{\rm (ii)} A necessary condition for 
$ \sum_{k=1}^{n} a_k\,\sin(a_kx) $ to be NN in some
neighborhood $ [0,\epsilon ] $, $ 0<\epsilon <\pi  $ is 
\begin{equation}  \sum_{k=1}^{n} ka_k \geq 0  \Label{bel2}  \end{equation}
\end{LEMMA}
\begin{PROOF}
Let us prove (i).
By assumption
\begin{equation}  0 \leq   \sum_{k=1}^{n} \frac{a_k\,\sin(kx)}{\pi -x}  \end{equation}
for all $ x\in[\pi -\epsilon ,\pi ) $
By taking the limit as $ x\rightarrow \pi  $, we get (using, for example, L'H\^opital's rule)
\begin{equation}  0 \leq   \lim_{x\rightarrow \pi }   \sum_{k=1}^{n} \frac{a_k\,\sin(kx)}{\pi -x} = \sum_{k=1}^{n} (-1)^{k+1}a_k .  \end{equation}

The proof of (ii) is similar.
\end{PROOF}

\begin{REM} \em
In the hypotheses of the Lemma, $ a_k $ are not required to be of the same
sign or monotone. Also note that unlike in the Belov condition, we are 
assuming in the hypothesis only that the sine polynomial itself 
(not any of its proper partial sums) is NN, and only one inequality
(\ref{bel1}) is required to hold (not for all $ n $).
\end{REM}

\begin{REM} \em
As Belov already pointed out, his condition (\ref{bel}) is no longer sufficient without 
the additional monotonicity requirement on the coefficients. We give an example
related to our sum $ \Phi  $. It is easy to verify that the polynomial
$$  2\sin(x)+\sin(2x)+ \frac{4}{3} \sin(3x) + \sin(4x) + \frac{6}{5} \sin(5x) + \frac{6}{8} \sin(8x)  $$
is not NN in $ [0,\pi ] $, although it satisfies (\ref{bel}).
This polynomial is constructed by taking $ \mb\Phi (5) $,
the first five terms of $ \mb\Phi  $, skipping
the terms involving $ \sin(6x) $ and $ \sin(7x) $ and add the next term with
a suitable coefficient to satisfy (\ref{bel}). The same is true for the polynomial constructed using
$ \mb\Phi (5) $ and $ \sin(10x) $. However, we notice that,
after that, all polynomials of the form
$$  \mb\Phi (5)+ \frac{6}{n} \,\sin(nx), \quad  n=12,14,16,\cdots  $$
are all PS.
\end{REM}

\begin{REM} \em
Another natural question to ask is whether
our Theorem~1 has a cosine counterpart, namely, whether $ \sum\,{\gamma _k}\,\cos(kx) $
is a PS, if $ \gamma _0=\gamma _1 $ and $ \gamma _k $ is given by (\ref{ak1}) for $ k=1,2,\cdots $. The answer
is also no. For $ x=\pi  $, the cosine series becomes
$$  \gamma _2 - \gamma _3 + \gamma _4 - \gamma _5 + \cdots  $$
and every partial sum with an even number of terms is negative, because $ \gamma _2<\gamma _3 $,
$ \gamma _4<\gamma _5 $, etc. A similar observation applies to the analogous sum 
$ \sum\,{\delta _k}\,\cos(kx) $.

\end{REM}

{\bf Acknowledgments}
The author is thankful to Horst Alzer for many inspiring discussions on
the subject of inequalities, in particular, trigonometric inequalities. 

Many of the technical
computations mentioned in the article were carried out using the
excellent Maple symbolic computation software.


\vspace{0.7cm}

\begin{thebibliography}{99}
% \bibitem{AAR} G.E. Andrews, R. Askey, R. Roy, Special Functions, Camb. Univ.
% Press, Cambridge, 1999.
\Bibitem{akl}
H. Alzer, S. Koumandos, and M. Lamprecht, M, A refinement of Vietoris'
inequality for sine polynomials, Math. Nachr. 283 (2010), 1549--1557.

\Bibitem{AK} H. Alzer, Man Kam Kwong, Sturm theorem and a reﬁnement of 
 Vietoris’ inequality for cosine polynomials,
arXiv:1406.0689 (math.CA).

\Bibitem{AK2} H. Alzer, Man Kam Kwong,
Rogosinski-Szeg\"o type inequalities for trigonometric sums,
J. Approx. Theory 190 (2015), 62--72.
% . ad
% R. Andreani and D.K. Dimitrov, An extremal nonnegative sine polynomial, Rocky Mt. J. Math.
% 33 (2003), 759¿774.

% \Bibitem{AK1}
% H. Alzer, S. Koumandos, Remarks on  a sine polynomial, Arch. Math. 93 (2009),
% 475-479.

% \Bibitem{AK2}
% H. Alzer, S. Koumandos, Sharp estimates for various trigonometric sums,
% Analysis 33 (2012), 9-26.

% \Bibitem{AY}
% H. Alzer, Q. Yin, On trigonometric sums in two variables, Jaen J. Approx. 4 (2012), 157-170.

\Bibitem{A} R. Askey, Orthogonal Polynomials and Special Functions, 
Reg. Conf. Ser. Appl. Math. (vol. 21), SIAM, Philadelphia, PA, 1975.

\Bibitem{AG} R. Askey, G. Gasper, Inequalities for polynomials, in: The Bieberbach conjecture (A. Baernstein II, D. Drasin, P. Duren, A. Marden, eds.), Math. surveys and monographs (no 21), Amer. Math. Soc., Providence, RI, 1986, pp. 7-32.

\Bibitem{AS} Askey, R., and Steinig, J.,
Some positive trigonometric sums,
Trans. Amer. Math. Soc. 187 (1074), 295-307.

\Bibitem{be} 
Belov, A.S., Examples of trigonometric series with nonnegative partial sums, Math. USSR Sb. 186, 21¿46 (1995) (Russian); 186, 485¿510 (1995) (English translation).

\Bibitem{b1} Brown, G.,
Positivity and boundedness of trigonometric sums,
Analysis in Theory and Applications 23 (2007), 380-388.

% \Bibitem{BH}
% G. Brown, E. Hewitt, A class of positive trigonometric sums, Math. Ann. 268 (1984),
% 91-122.

% \Bibitem{DM}
% D.K. Dimitrov, C.A. Merlo, Nonnegative trigonometric polynomials, Constr. Approx. 18 (2002), 117-143.
% 
% . du
% B.A. Dumitrescu, {\em Positive Trigonometric Polynomials and Signal
% Processing Applications}, Springer, Berlin (2007)
% 
% \Bibitem{f} L. Fej\'er, Einige S\"atze, die sich auf das Vorzeichen einer 
% ganzen rationalen Funktion beziehen, Monatsh. Math. Phys. 35 (1928), 305-344.
% 
% . fd
% J.J. Fern\'andez-Dur\'an, Circular distributions based on nonnegative
% trigonometric sums, Bio-
% metrics 60 (2004), 499¿503.
% 
% \Bibitem{G}
% G. Gasper, Nonnegative sums of cosine, ultraspherical and Jacobi polynomials, J. Math. Anal. Appl. 26 (1969), 60-68.
% 
% . gh A. Gluchoff and F. Hartmann, Univalent polynomials and non-negative
% trigonometrics sums,
% Am. Math. Mon. 195 (1998), 508¿522 
% 
% \Bibitem{Ks}
% S. Koumandos, Some inequalities for cosine sums,
% Math. Inequal. Appl. 4 (2001), 267-279.

\Bibitem{Ks2}
S. Koumandos, Inequalities for trigonometric sums,
{\em Nonlinear Analysis},
Springer Optimization and Its Applications Volume 68, 2012, pp 387-416

\Bibitem{Kw}
Man Kam Kwong, Nonnegative trigonometric
polynomials, Sturm's theorem, and symbolic
computation, arXiv: 1402.6778 [math CA]
(2014).

\Bibitem{Kw2} Man Kam Kwong,
An improved Vietoris sine inequality. J. Approx. Theory 189 (2015), 29--42. 

% \bibitem{MMR} G.V. Milovanovi\'c, D.S. Mitrinovi\'c, Th.M. Rassias, Topics in Polynomials: Extremal Problems, Inequalities, Zeros, World Sci. Publ.,
% Singapore, 1994.

% \bibitem{S} G. Szeg\"o, Power series with multiply monotonic sequences of 
% coefficients, Duke Math. J. 8 (1941), 559-564.
% 
% \bibitem{T} P. Tur\'an, \"Uber die arithmetischen Mittel der Fourierreihe,
% J. London Math. Soc. 10 (1935), 277-280.

\Bibitem{V} L. Vietoris, \"Uber  das  Vorzeichen  gewisser  Trigonometrische  Summen,  S.-B. Osterreich.  Akad.  Wiss.,
167 (1958), 125-135. Teil II: Anz.
Osterreich. Akad. Wiss., (1959), 192-193.

\end{thebibliography}
\end{document}